\documentclass[10pt,twoside]{article}
\usepackage{graphicx}
\usepackage{amsmath}
\def\dref#1{(\ref#1)}

\def\ra{\longrightarrow}

\def\deq{\stackrel{\Delta}{=}}
\def\R{I\!\!R}
\def\F{{\cal F}}

\def\dlim{\displaystyle\lim\limits}
\def\dsup{\displaystyle\sup\limits}
\def\argmin{\mathop{\mbox{argmin}}}
\def\s{\theta}
\def\hatf{{\hat f}}

\def\F{{\cal F}}

\def\argmin{\mathop{\mbox{argmin}}}

\def\ra{\longrightarrow} \def\R{{\cal R}}

\def\dsum{\mathop{\sum}\limits} 
\def\dsup{\mathop{\sup}\limits}
\def\dlim{\mathop{\lim}\limits}

\def\deq{\stackrel{\Delta}{=}}
\def\underb{{\underline b}}
\def\overb{{\overline b}}
\def\hatf{{\hat f}}
\def\ba{\begin{array}}
\def\ea{\end{array}}
\def\be{\begin{equation}}

\def\ee{\end{equation}}
\usepackage{Latex-document}

\title{\bf Exploring the Capability and Limits\vskip -2mm of the Feedback Mechanism
\vskip 6mm}

\author{Lei Guo\vspace*{-0.5cm}\thanks{Institute of Systems Science, Chinese Academy of
Sciences, Beijing 100080, China. E-mail: Lguo@control.iss.ac.cn}}
\date{\vspace{-8mm}}

\begin{document}

\maketitle

\thispagestyle{first} \setcounter{page}{785}

\begin{abstract}

\vskip 3mm

Feedback is a most important concept in control systems, its main
purpose is to deal with internal and/or external uncertainties in
dynamical systems, by using the on-line observed information.
Thus, a fundamental problem in control theory is to understand the
maximum capability and potential limits of the feedback mechanism.
This paper gives a survey of some basic ideas and results
developed recently in this direction, for several typical classes
of uncertain dynamical systems including parametric and
nonparametric nonlinear systems, sampled-data systems and
time-varying stochastic systems.

\vskip 4.5mm

\end{abstract}

\vskip 12mm

\section{Introduction} \label{section 1}\setzero
\vskip-5mm \hspace{5mm }

Feedback is ubiquitous, and exists in almost all goal-directed
behaviors [1]. It is indispensable to the human intelligence, and
is important in learning, adaptation, organization and evolution,
etc. Feedback is also the most important concept in control, which
is a fundamental systems principal when dealing with uncertainties
in complex dynamical systems. The uncertainties of a system are
usually classified into two types: internal (structure) and
external (disturbance) uncertainties, depending on the specific
dynamical systems to be controlled. Feedback needs information,
and there are also two types of information in a control system:
{\it a priori} information and {\it posteriori} information. The
former is the available information before controlling a system,
while the later is the information exhibited by the system dynamic
behaviors. It is the {\it posteriori} information that makes it
possible for the feedback to reduce the influences of the
uncertainties on control systems. Two of the fundamental questions
in control theory are: How much uncertainty can be dealt with by
feedback ? What are the limits of feedback ? These are conundrums,
despite of the considerable progress in control theory over the
past several decades.

    The existing feedback theory in control systems can be
roughly classified into three groups: traditional feedback, robust
feedback and adaptive feedback. In the ideal case where the
mathematical model can exactly describe the true system, the
feedback law that are designed based on the full knowledge of the
model may be referred to as traditional feedback. Unfortunately,
as is well known, almost all mathematical models are
approximations of practical systems, and in many cases there are
inevitable large uncertainties in our mathematical descriptions.
The primary motivation of robust and adaptive control is to deal
with uncertainties by designing feedback laws, and much progress
has been made in these two areas. Robust feedback design allows
that the true system model is not exactly known but lies in a
``ball" centered at a known nominal model with reliable model
error bounds(cf. e.g. \cite{2}, \cite{4}).

By adaptive feedback we mean the (nonlinear) feedback which
captures the uncertain (structure or parameter) information of the
underlying system by properly utilizing the measured on-line data.
The well-known certainty-equivalence principle in adaptive control
is an example of such philosophy. Since an on-line learning
mechanism is usually embedded in the structure of adaptive
feedback, it is conceivable that adaptive feedback can deal with
larger uncertainties than other forms of feedback can do. Over the
past several decades, much progress has been made in the area of
adaptive control (cf. e.g. \cite{6}--\cite{11}). For linear finite
dimensional systems with uncertain parameters, a well-developed
theory of adaptive control exists today, both for stochastic
systems (cf. \cite{7},\cite{10},\cite{11}) and for deterministic
systems with small unmodelled dynamics (cf. \cite{8}). This theory
can be generalized to nonlinear systems with linear unknown
parameters and with linearly growing nonlinearities (cf.
e.g.\cite{16}). However, fundamental difficulties may emerge in
the design of stabilizing adaptive controls when these structural
conditions are removed. This has motivated  a series of studies on
the maximum capability (and limits) of the feedback mechanism
starting  from \cite{14}.

To explore the maximum capability and potential limits of the
feedback mechanism, we have to place ourselves in a framework that
is somewhat different from the traditional robust control and
adaptive control. First, the system structure uncertainty may be
nonlinear and/or nonparametric, and a known or reliable ball
containing the true system, which is centered at a known nominal
model, may not be available {\it a priori}. Second, we need to
consider the maximum capability of the whole feedback mechanism
(not only a fixed feedback law or a special class of feedback
laws). Moreover, we need to answer not only  what the feedback can
do, but also the more difficult and important question, what the
feedback can not do. We shall also work with discrete-time (or
sampled-data) feedback laws, as they can reflect the basic causal
law as well as the limitations of actuator and sensor in a certain
sense, when implemented with digital computers.

In this talk, we will give a survey of some basic ideas and
results obtained in the recent few years (\cite{14}--\cite{21}),
towards understanding the capability and limits of the feedback
mechanism in dealing with uncertainties. For several basic classes
of uncertain nonlinear dynamical control systems, we will give
some critical values and ``Impossibility Theorems" concerning the
capability of feedback. The reminder of the paper is organized as
follows: the problem formulation will be given in Section 2, and
some basic classes of discrete-time parametric and nonparametric
nonlinear control systems will be studied in Sections 3 and 4,
respectively. Other classes of uncertain systems, including
sampled-data systems and time-varying linear systems with hidden
Markovian jumps, will be considered in Section 5, and some open
problems will be stated in the concluding remarks in Section 6.

\section{Problem formulation} \label{section 2}
\setzero\vskip-5mm \hspace{5mm }

Let $\{u_t, t\geq 0\}$ be an ${\cal R}^m$-valued input process of
a discrete-time or sampled-data uncertain dynamical  control
system (whose structure is unknown or not fully known, and is
subject to some noise disturbances), and let $\{y_t, t\geq 0\}$ be
the corresponding on-line observed $\R^n$-valued output process
(see the following figure).

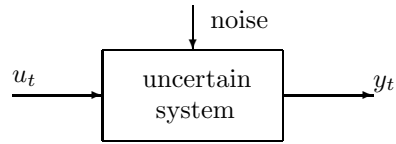
\centerline{ \setlength{\unitlength}{0.6cm}
\begin{picture}(8,4)
\put(0,1){\vector(1,0){2}} \put(0,1.3){\mbox{$u_t$}}
%\put(2,1){\line(1,0){2}}
\put(2,2){\line(0,-1){2}} \put(2,2){\line(1,0){4}}
\put(6,2){\line(0,-1){2}} \put(2,0){\line(1,0){4}}
\put(6,1){\vector(1,0){2}} \put(8,1.2){\mbox{$y_t$}}
\put(4,3){\vector(0,-1){1}} \put(4.4,2.5){\mbox{noise}}
 \put(2.5,1.2){\mbox{
\,\,uncertain}} \put(2.5,0.5){\,\quad system}
\end{picture}
}

\vskip 0.2cm
 At any moment $t \geq 0$, the input signal $u_t$ is
said to be a \it feedback signal, \rm if there is a Lebesgue
measurable mapping $h_t(\cdot): \R^{(t+1)n} \to \R^m $ such that
$$
u_t = h_t(y_0, y_1, ..., y_t).
$$
In other words, $u_t$ is a function of the {\it posteriori}
information observed up to time $t$. A \it feedback law \rm $u$ is
a sequence of feedback signals, i.e., $u = \{u_t, t\geq 0\}$.
Furthermore, the \it feedback mechanism \rm $U$ is defined as the
set of all possible feedback laws,
$$
U = \{u |\,\,\, u \quad \mbox{is any feedback law} \}.
$$

For a complex system whose structure contains uncertainties, it is
not a simple (and in fact difficult) problem to find a feedback
law $u$, such that the corresponding output process can achieve a
desired goal. This involves questions like: what kind of
properties of an uncertain dynamical system can be changed by
feedback ? how to construct a satisfactory feedback based on the
available information ? More fundamental questions are: how much
uncertainty can be deal with by feedback ?  What are the maximum
capability and limits of the feedback mechanism $U$ ? In the next
three sections, we will present some preliminary results towards
answering these questions. For further discussion, we need the
following definition.

{\bf Definition 2.1}. {\it A dynamical control system is said to
be globally stabilizable if there exists a feedback law $u \in U$,
such that the output process of the system is bounded in the mean
square sense, i.e.,
$$
\sup_{t \geq 0} E\|y_t\|^2<\infty, \quad \mbox{ for any initial
value}\,\,y_0.
$$}

\section{Parametric nonlinear systems} \label{section 3}
\setzero\vskip-5mm \hspace{5mm }

Consider the following basic discrete-time parametric nonlinear
system:
\begin{align}\label{3.1}
  y_{t+1} = \theta f(y_t) + u_t + w_{t+1},
\end{align}
where $y_t$, $u_t$ and  $w_t$ are the scalar  system output, input
and noise processes, respectively. For simplicity, we assume that

A1) $\{w_t\}$ is a Gaussian noise process;

A2) $\theta$ is an unknown non-degenerate Gaussian random
parameter;

A3)  The function $f(\cdot)$ is known and has the following growth
rate:
$$
f(x) \quad \sim \quad M x^b,  \quad\mbox{as}\quad x\to\infty,
$$
where $b\geq 0$, $M>0$ are constants. Obviously, if $b\leq 1$,
then the nonlinear function $f(\cdot)$ has a growth rate which is
bounded by linear growth. This case can be easily dealt with by
the existing theory in adaptive control (see e.g. \cite{16}). Our
prime concern here is to know whether or not the system can be
globally stabilized by  feedback for any $b>1$?

The following theorem gives a critical value of $b$, which
characterizes the maximum capability of the feedback mechanism.

{\bf Theorem 3.1.}{\it Consider the system (\ref{3.1}) with
Assumptions A1)--A3) holding. Then $ b=4$  is a critical case  for
feedback stabilizability. In other words,

  (i). If $b\geq 4$, then for
any feedback law  $u \in U$, there always exists a set $D$ (in the
basic probability space) with positive probability such that
$$
|y_t| \to\infty, \mbox{ on $D$}
$$
at a rate faster than exponential.

 (ii). If $b < 4$, then the least-squares-based adaptive minimum
 variance feedback
control $u_t=-\theta_t f(y_t)$ where $\theta_t$ is the
least-squares estimate for $\theta$ at time $t$, can render the
system to be globally stable and optimal, with the best rate of
convergence:
$$
\sum^T_{t=1} (y_t  - w_t)^2 = O(\log T), \mbox{ a.s., as } T\ra
\infty.
$$}
\vskip 0.1cm

\bf Remark 3.1. \rm This result is somewhat surprising since the
assumptions in our problem formulation  have no explicit
relationships with the value  $b=4$. We remark that the related
results were first found and established in a somewhat general
framework in \cite{14}. In particular, the first part (i) was
contained in Remark~2.2 in \cite{14}, and was later extended to
general unknown parameter case in \cite{16} by using a conditional
Cramer-Rao inequality. The second part (ii) is a special case of
Theorem~2.2 in \cite{14}.

\bf Remark 3.2. \rm There are many implications of Theorem~3.1.
For example, the limitation of feedback given in Theorem 3.1 (i)
is readily applicable to general class of  uncertain systems of
the form
$$
y_{t+1} = f_t(y_t, ..., y_{t-p}, u_t, ..., u_{t-q}) + w_{t+1},
$$
as long as it contains the basic class (\ref{3.1}) as a subclass.
Theorem 2.1 can also be used to show the fundamental differences
between continuous-time and discrete-time nonlinear adaptive
control(see also, \cite{16}). Note that the noise free case can be
trivially controlled, regardless of the growth rate of the
nonlinear function $f(\cdot)$ (see \cite{15}). This means that the
noise effect in (\ref{3.1}) plays an essential role in the
non-stabilizability result of Theorem 3.1 (i): the noise effect
gives estimation errors to even the ``best'' parameter estimates,
which are then amplified step by step by the nonlinearity of the
system, leading to the final instability of the closed-loop
systems, despite of the strong consistency of the parameter
estimates (\cite{14}) .

 Theorem 3.1 concerns with the case where the unknown
parameter $\theta $ is a scalar. To see what happens when the
number of the unknown parameters increases, let us consider the
following polynomial nonlinear regression:
\begin{align}\label{2}
y_{t+1}=\theta_1y_t^{b_1}+\theta_2y_t^{b_2}+\cdots +\theta_py_{t}^{b_p}+u_t+w_{t+1}.
\end{align}
Again, for simplicity, we assume that

A1)$'$ $b_1>b_2>\cdots>b_p>0$;

A2)$'$ $\{w_t\}$ is a sequence of independent random variable with
a common distribution $N(0,1)$;

A3)$'$ $\theta\deq [\theta_1\cdots\theta_p]^\tau$ is   a random
parameter with distribution  $  N({\overline\theta}, I_p)$.

Now, introduce a characteristic polynomial
$$
  P(z)=z^{p+1}-b_1 z^p+(b_1-b_2)
z^{p-1}+\cdots+(b_{p-1}-b_p)z+b_p,
$$
which plays a crucial role in characterizing the limits of the
feedback mechanism as shown by the following ``impossibility
theorem''.

{\bf Theorem 3.2.} {\it If there exists a real number $z\in (1,b_1)$ such that $ P(z)<0$, then the above system
\dref{2} is not stabilizable by feedback. In fact, for any feedback law $ u \in U$ and any initial condition
$y_0\in \R^1$, it is always true that
$$
  E|y_t|^2\rightarrow \infty, \mbox{ as } t\rightarrow \infty
 $$
at a rate faster than exponential.}

\bf Remark 3.3. \rm The proof of Theorem 3.2 can be found in
\cite{15}, and extensions to non-Gaussian parameter case can be
found in \cite{16}. An important consequence of this theorem is
that the system \dref{2} is not stabilizable by feedback in
general, whenever $b_1
> 1$ and the number of unknown parameters $p$ is large(see \cite{15}).
This fact implies that the uncertain nonlinear system
$$
y_{t+1}=f(y_t)+u_t+w_t
$$
with  $f(\cdot)$ being unknown but satisfying
$$
\|f(x)\|\leq c_1+c_2\|x\|^b, \quad b>1,
$$
may not be stabilizable by  feedback in general. This gives us
another fundamental limits on feedback in the presence of
parametric uncertainties in nonlinear systems, and motivates the
study of nonparametric control systems with linear growth
conditions in the next section.

\section{Nonparametric nonlinear systems} \label{section 4} \setzero\vskip-5mm \hspace{5mm }

Consider the following first-order nonparametric   control system:
\begin{align}\label{3}
y_{t+1}=f(y_t) +u_t +w_{t+1}, \quad t\geq 0,\quad y_0 \in \R^1 ,
\end{align}
where   $\{y_t\}$ and $\{u_t\}$ are the scalar output and input,
and $\{w_t\}$ is  an ``unknown but bounded'' noise  sequence, i.e.
$|w_t|\leq w$, $\forall t$, for some constant $w>0$. The nonlinear
function $f(\cdot): \,\R^{1} \rightarrow \R^{1}$ is assumed to be
completely unknown. We are interested in understanding
 how much uncertainty in $ f(\cdot) $  can
be dealt with by feedback. In order to do so, we need to introduce
a proper measure of uncertainty first.

Now, define $ \F \deq \{f: \, \R^{1} \rightarrow \R^{1} \} $ and
introduce a quasi-norm on $\F$ as follows:
$$
\|f\| \deq \dlim_{\alpha \rightarrow \infty} \dsup_{(x,y)\in \R^2}
     \dfrac{|f(x)-f(y)|}{|x-y|+\alpha}, \quad \forall f \in \F.
$$

Having introduced the norm $\|\cdot\|$, we can then define a ball
in the space $(\F, \|\cdot\|)$ centered at its ``zero'' $\theta_F$
with radius $L$:
$$
\F(L) \deq \{f \in \F:\, \|f \| \leq L \}
$$
where $\theta_F \deq \{f \in \F: \, \|f\|=0 \}$. It is obvious
that the size of $\F(L)$ depends on the radius $L$, which will be
regarded as the measure of the size of uncertainty in our study to
follow.

The following theorem establishes a quantitative relationship
between the capability of feedback and the size of uncertainty.

{\bf Theorem 4.1.}{\it Consider the nonparametric control system
\dref{3}. Then
  the maximum uncertainty  that can be dealt with by feedback is
a ball with radius $L = \frac32 +\sqrt{2}$ in the normed function
space $(\F, \|\cdot\|)$, centered at the zero $\theta_F$. To be
precise,

 (i) If $L<\frac32+\sqrt{2}$, then there exists a feedback law
 $u\in U$
  such that for any $f\in \F(L)$, the corresponding
 closed-loop control system \dref{3} is globally stable in the
 sense that
 $$
 \sup_{t\geq 0} \{|y_t|+|u_t|\}<\infty, \quad \forall y_0\in \R^1;
 $$

(ii) If $L\geq \frac32+\sqrt{2}$, then for any  feedback law $u\in
U$
  and any initial value $y_0\in \R^1$, there always exists some $f\in
\F(L)$ such that the corresponding closed-loop system \dref{3} is
unstable, i.e.,
$$
\sup_{t\geq 0} |y_t|=\infty.
$$}

The proof of the above theorem is given in \cite{18}, where it is
also shown that once the stability of the closed-loop
  system is established, it is a relatively easy task to
evaluate the control performance.

\bf Remark 4.1. \rm The stabilizing feedback law  in Theorem~4.1
(i)
  can be constructed as follows(see \cite{18}):

$$
u_t=\left\{
\begin{array}{ll}
u'_t, \qquad & \mbox{if } |y_t-y_{i_t}|>\epsilon\\
u''_t, & \mbox{if } |y_t-y_{i_t}|\leq \epsilon
\end{array}
\right.
$$
where $\epsilon>0$ is any given threshold. In other words, $u_t$
is a switching feedback based on a
 stabilizing feedback $u'_t$ and a tracking feedback $u''_t$,
 which are defined as follows:
     $$u'_t=-\hatf_t(y_t)+\frac12(\underb_t+\overb_t)$$
     where $\hatf$ is  the
  nearest neighbor (NN) estimate of $f$ defined by
$$
{\hat f}_t(y_t)\deq y_{i_t+1}-u_{i_t}
$$
 with
$$
i_t\deq \argmin_{0\leq i\leq t-1} |y_t-y_i|
$$
and where
   $$
\underb_t=\min_{0\leq i\leq t} y_i, \quad \overb_t=\max_{0\leq
i\leq t} y_i
$$

The tracking feedback $u_t''$ is defined
 $$u_t''=-\hatf_t(y_t)+y^*_{t+1},$$
 where $\{y^*_t\}$ is a bounded reference sequence.
It is obvious that $u_t$ depends on the observations $\{y_0,
y_1,\ldots,y_t\}$ only.

\vskip 0.3cm
 One may try to generalize
Theorem 4.1 to the following high-order nonlinear systems ($p\geq
1$):
\begin{align}\label{4}
y_{t+1}=f(y_t, y_{t-1}, \ldots, y_{t-p+1})+u_t+w_{t+1}
\end{align}
where $f(\cdot): R^p\ra R^1$ is assumed to be completely unknown,
but belongs to the following class of Lipschitz functions:
$$
\F(L)=\{f(\cdot):|f(x)-f(y)|\leq L\|x-y\|,\forall x, y\in \R^p\}
$$
where $L>0$, $\|x\|=\dsum_{i=1}^p|x_i|$,
$x=(a_1,\ldots,x_p)^\tau\in R^p$. Again, $\{w_t\}$ is a sequence
of ``unknown but bounded'' noises. The following ``impossibility
theorem" is established in \cite{21}. \vskip 0.3cm

{ \bf Theorem 4.2. }{\it If $L$ and $p$ satisfy
\begin{align}\label{5}
L+\frac12\geq (1+\frac1p)(pL)^{\frac1{p+1}}
\end{align}
then there does not exist any globally stabilizing feedback law
for the class of uncertain systems \dref{4} with  $f\in \F(L)$.}

It is easy to see that if $p=1$ then the above inequality \dref{5}
reduces to  $L\geq \frac32+\sqrt{2}$, which we know to be a
critical value for the capability of feedback by Theorem~4.1.
However, when $p>1$ and \dref{5} does not hold, whether or not
there exists a stabilizing   feedback for the uncertain system
\dref{4} with $f\in \F(L)$ still remains as an open question.

\section{Other uncertain systems  } \label{section 5}
\setzero\vskip-5mm \hspace{5mm }

 In this section, we briefly mention some related results on other basic classes
 of uncertain systems.

 Let us first consider the following simple but basic continuous-time system:
   \begin{equation}
   \label{6}
    {\dot x_{t}}=f(x_t)+u_t, \quad t\geq 0, x_0 \in \R^1.
   \end{equation}

The system signals are assumed to be sampled at a constant rate
$h>0$, and the input is assumed to be implemented via the familiar
zero-order hold device (i.e., piecewise constant functions):
   \begin{equation}
   \label{7}
u_t=u_{kh},\quad kh\leq t<(k+1)h,
   \end{equation}
where $u_{kh}$ depends on $\{ x_0,x_h,...,x_{kh}\}$.

The nonlinear function $f$ in (\ref{6}) is assumed to be unknown
but belongs to the following class of local Lipschitz (LL)
functions:
   \begin{equation}
   \label{8}
{G^L_c}=\{ f | \,f\mbox{ is  LL and satifies}\\
       \;    |f(x)|\leq L|x|+c, \forall x \in R^1 \},
   \end{equation}
where $c>0$ and $L>0$ are constants.  The ``slope'' $L$  of the
unknown nonlinear functions in $G^L_c$ may be regarded as a
measure of the size of the uncertainty. Similar to the
discrete-time case in Theorem~4.1, $L$ plays a crucial role in the
determination of  the capability and limits of the sample-data
feedback \cite{17}.

  {\bf Theorem 5.1.} {\it Consider the sampled-data control system
  \dref{6}--\dref{7}.
 If $Lh>7.53$, then for any $c>0$ and any sampled-data control $\{
u_{kh},k\geq 0\}$ there always exists a function $f^* \in G^L_c$,
such that the state signal of (\ref{6})--(\ref{7})
 corresponding to $f^*$ with initial point $x_0=0$ satisfies $(k \geq 1)$ }
        $$|x_{kh}|\geq (\frac {Lh}{2})^{k-1} \cdot ch
              \mathop{\longrightarrow}\limits_{ k\rightarrow\infty } \infty.$$

 \bf Remark 5.1. \rm  This ``impossibility theorem'' shows that if $Lh$
is larger than a certain value, then there will exist no
stabilizing sampled-data feedback. On the other hand, it is easy
to show that if $Lh<\log 4$, then a globally stabilizing
sampled-data feedback can be constructed (see \cite{17}). An
obvious open question here is how to bridge the gap between $\log
4$ and 7.53. Needless to say, Theorem 5.1 gives us some useful
quantitative guidelines in choosing properly the sampling rate in
practical applications. \vskip 0.1cm

Next, we consider the following linear time-varying stochastic
model:
\begin{align}
\label{101}
 \ba{l}
x_{t+1}=A(\s_t)x_t +B(\s_t)u_t +w_{t+1}, \, \, t\geq 1; \\
\ea
\end{align}
where $x_t \in \R^n$ ,$u_t\in \R^m$ and $w_{t+1}\in \R^n $ are the
state, input and noise vectors respectively.  We assume that

\vskip 1mm H1) $\{ \s_t \}$ is an unobservable Markov chain which
is homogeneous, irreducible and aperiodic, and which takes values
in a finite set $\{ 1,2,\cdots,N \}$ with transition matrix
denoted by $P=(p_{ij})_{N\*N}$, where by definition $p_{ij}=
P\{\s_t =j|\s_{t-1} =i\}$.

\vskip 1mm H2) There exists some $m \times n$ matrix $K$ such that
$ det \big[ (A_i-A_j)-(B_i-B_j)K \big] \not= 0, \quad \forall
i\not= j, 1 \leq i,j \leq N$, where $A_i\deq A(i) \in \R^{n\* n}$,
$B_i\deq B(i) \in \R^{n\* m}$ are the system matrices.

\vskip 1mm H3) $\{w_t\}$ is a martingale difference sequence which
is independent of $\{ \s_t\}$, and satisfies $
 {\sigma I \leq E w_t w_t', \,\quad Ew_t'w_t \leq
\sigma_w,\quad \forall t} $,  where $\sigma$ and $\sigma_w$ are
two positive constants, and the prime superscript represents
matrix transpose. \vskip 1mm

For simplicity of presentation, we denote $S\deq
\{1,2,\cdots,N\}$.  The following theorem gives a  fairly complete
characterization of feedback stabilizability for the hidden
Markovian model (\ref{101})\cite{20}.

\vskip 3mm {\bf Theorem 5.2.} {\it Let the above Assumptions
H1)--H3) hold for the dynamical system (\ref{101}) with hidden
Markovian switching. Then the system is stabilizable by feedback
{\it if and only if}   the following coupled algebraic
Riccati-like equations have a solution consisting of $N$ positive
definite matrices $\{ M_i >0, i \in S \}$:
$$
   \dsum_{j}A_j' p_{ij}M_j A_j  -
    \big( \dsum_{j}A_j' p_{ij}M_j B_j \big)
    \big( \dsum_{j}B_j' p_{ij}M_j B_j \big)^+
    \big( \dsum_{j}B_j' p_{ij}M_j A_j \big)
                    -M_i=-I,
$$
 where $i\in S$ and $(\cdot)^+$ denotes the Moore-Penrose
generalized-inverse of the corresponding matrix.} \vskip 0.2cm

\bf Remark 5.2. \rm Theorem 5.2 shows that the capability of
feedback depends on both the structure complexity measured by
$\{A_j, B_j, 1\leq j\leq N\}$ and the information uncertainty
measured by $\{p_{ij}, 1\leq i, j\leq N\}$. To make it more clear
in understanding how the capability of feedback depends on both
the complexity and uncertainty of the system, we consider the
simple scalar variable case with $B(\theta_t)=1$, where the Markov
chain $\{\theta_t\}$ has two states $\{1,2\}$ only and
$p_{12}=p_{21}$. It can be shown by Theorem 5.2 that the system is
stabilizable if and only if $CP<1$, where $C\deq (A_2-A_1)^2$ and
$P\deq (1-p_{12})p_{12}$ can be interpreted as measures of the
structure complexity (degree of dispersion) and the information
uncertainty respectively (see \cite{19} for details).

\section{Concluding remarks}\label{section 6} \setzero\vskip-5mm \hspace{5mm }

For several basic classes of uncertain dynamical control systems,
we have given some critical values or equations to characterize
the capability and limits of the feedback mechanism, and have
shown that``impossibility theorems'' hold even for some seemingly
simple uncertain dynamical systems. Of course, many important
problems still remain open. Examples are as follows:

(i)  For general high-dimensional or high-order uncertain
nonlinear control systems, to find critical conditions
characterizing the capability of feedback, at least to find
general sufficient conditions under which feedback stabilization
is possible in the discrete-time case.

(ii) To characterize the maximum capability of feedback that is
designed based on switched linear control models, in dealing with
uncertain nonlinear dynamical systems.

(iii) To find a suitable mathematical framework within which the
issue of establishing a quantitative relationship among {\it a
priori} information, feedback performance and computational
complexity can be addressed adequately.

\label{lastpage}

\end{document}